\documentclass{article}
\usepackage{amsmath}
\usepackage{amsfonts}
\usepackage{amssymb}
\usepackage{euscript}
\usepackage{epsfig}

\def\dj{d\kern-0.4em\char"16\kern-0.1em}

 \newtheorem{thm}{Theorem}[section]
 
 \newtheorem{lem}[thm]{Lemma}
 \newtheorem{cor}[thm]{Corollary}

 \newtheorem{rem}{Remark}[section]

\title{
POWER SERIES DETERMINED BY \\
AN EXPERIMENT ON THE UNIT INTERVAL
}
\author{
{\bf Predrag M. Rajkovi\'c$^1$\footnote{The corresponding author, {\it E-mail}:\ pedja.rajk@masfak.ni.ac.rs}},
{\bf Sladjana D. Marinkovi\'c$^2$}\\
{\bf Miomir S. Stankovi\'c$^3$}\\
{\small $^1$Department of Mathematics, Faculty of Mechanical Engineering}\\
{\small $^2$Department of Mathematics, Faculty of Electronic Engineering}\\
{\small $^3$Department of Mathematics, Faculty of Occupational Safety}\\
[2mm] {\bf   University of Ni\v s,\ Serbia}
 }

\begin{document}
\maketitle

\begin{abstract}
We consider the linear combinations of elements of two sequences: the first one a priory given nonnegative sequence and the second random sequence from the unit interval.  We investigate the expected value of  the  smallest natural number such that  the value of these linear combinations exceed a positive number. After very clear geometrical conclusions, we find the function which expresses the expected value. Here, we recognize a few known results like the special cases.  \end{abstract}

\smallskip

{\bf MSC 2010}: 34K60, 60G50, 33B10.

{\bf Keywords:} Experiment, expected value, probability, exponential function.


\section{Introduction}
The well known number $e\approx 2.71818\ldots$ appeared like an expected value in a trial problem by  Putnam \cite{Bush} in 1958. This wonderful, riddle-like problem motivated others to solve it in various ways \cite{Shultz}. Even more, B. \'Curgus and  R.I. Jewetts \cite{Curgus} have considered like a function of upper bound  which should be exceed. In that manner, they have got the exponential function.

	 In this paper, we will expose our generalization of this problem.
	 Namely, consider the next experiment.

\smallskip

{\bf Procedure.}

Let $ a = \{a_k\}_{k\in\mathbb N}$ be a given positive nondecreasing sequence, i.e.
\begin{equation}\label{seq}
0<a_1\le a_2\le a_3\le \cdots \le a_k\le a_{k+1}\le \cdots \ .
\end{equation}
To any random sequence
$x=\{x_n\}\ (x_n\in[0,1], \forall n\in\mathbb N)$, we can join the finite sequence of weighted partial sums
$$
S_{k,m}(a,x) = a_k x_k + a_2 x_2+\cdots + a_m x_m.
$$
Take $t>0$.

Step 1. Take a random number $x_1$ from the interval $[0,1]$ and denote $n=1$.

Step 2. If $S_{1,n}(x,a)>t$, then memorize the value $n$ and stop.

Step 3. If $S_{1,n}(x,a)\le t$, then increase $n$  into  $n+1$, take the next random number $x_n$ from $[0,1]$ and return to Step 2.

Repeat this procedure enough times. What is the expected  value for $n$?

We are looking for the expectation
\begin{equation}\label{function}
f(t, a)  =  E\{n\in\mathbb N:  S_{1,n-1}(x,a)\le t < S_{1,n}(x,a)\}.
\end{equation}
Denote with $p_m(t)$ the following function
\begin{equation}\label{pm}
p_m=\mu(P_m), \quad \text{\rm where}\quad P_m=\{x\in[0,1]^m:S_{1,m}(a,x)\le t\}\,.
\end{equation}
Here, $\mu$ is the Euclidean measure in $\mathbb R^m$.

The expected value for $n$ is
\begin{equation*}
f(t, a) = \sum_{m=0}^\infty m(p_{m-1}-p_m).
\end{equation*}
If it is fulfilled
\begin{equation}\label{req2}
\lim_{m\to \infty} mp_m =0,
\end{equation}
then
\begin{equation}\label{formula2}
f(t, a)=\sum_{n=0}^\infty p_n.
\end{equation}

Also, we will use notation
\begin{equation*}
(t)_+ =
\begin{cases}
t,& t > 0,  \\
   0,        &t \le 0,
\end{cases}
\qquad
f_+(t)=
\begin{cases}
f(t),& f(t) > 0,  \\
   0,        &f(t) \le 0.
\end{cases}
\end{equation*}

\section{The geometrical approach}

Here, we will expose cases which can be illustrated by one, two and three dimensional figures.

We accept $p_0=1$.

In one dimensional space, there are 2 cases:

\noindent(1.1) For $0\le t\le a_1$, i.e., when is valid $t/a_1\le 1$, the border point is $X_1(t/a_1)$. Hence $P_1=X_1O$, and the measure is
$$
p_1 = \mu(P_1) = \frac{t}{a_1}\ .
$$
(1.2) For $a_1<t$, it is $P_1 = [0,1]$, and $p_1=1$.

\smallskip

In two dimensional space, the line $a_1x_1+a_2x_2=t$ passes through the points $X_1(t/a_1,0),\ X_2(0,t/a_2)$.

\noindent(2.1) For $0\le t\le a_1$, the whole segment $X_2X_1$ belongs to the square $[0,1]^2$. Hence $P_2$ is the interior of the triangle $X_2X_1O$ and the measure is
$$
p_2 = \mu(X_2X_1O) = \int_0^{t/a_1} dx_1 \int_0^{(t-a_1x_1)/a_2} dx_2 = \frac{1}{2} \frac{t}{a_1} \frac{t}{a_2}\ .
$$

\noindent(2.2) Let $a_1< t\le a_2$.  This condition guarantees that the points $O(0,0)$ and $S(1,1)$ are from the different sides of the line $a_1x_1+a_2x_2=t$. That is why it exists the cutting segment of this line and the unit square. Since $t-a_1>0$, or equivalently $1<t/a_1$, the point $X_1$ is outside the unit square, like on the Figure 3a.  Hence
$$
p_2=\frac{t^2-(t-a_1)^2}{2a_1a_2}\,.
$$

\noindent(2.3) Let $a_2< t\le a_1+a_2$.
In that case (see Figure 3b), the segment $X_1X_2$ contains the points $B(1,(t-a_1)/a_2)$ and $C((t-a_2)/a_1,1)$. This line divides the square $[0,1]^2$ into 2 pieces. Hence
$$
P_2 = \overline{OS_{1,0}BCS_{0,1}} = \overline{X_2X_1O}\setminus\bigl( \overline{X_1S_{1,0}B}\cup \overline{X_2CS_{0,1}}\bigr),
$$
wherefrom
$$
p_2=\mu(P_2)=\frac12\frac{t}{a_1}\frac{t}{a_2}
-\frac12\Bigl(\frac{t}{a_1}-1\Bigr)\frac{t-a_1}{a_2}-\frac12\Bigl(\frac{t}{a_2}-1\Bigr)\frac{t-a_2}{a_1},
$$
i.e.
$$
p_2=\frac{t^2-(t-a_1)^2-(t-a_2)^2}{2a_1a_2}\,.
$$

\smallskip

We can unify (2.1-3) into one formula:
$$
p_2(t) = \frac1{2a_1a_2} \Bigl(t^2 - \sum_{k=1}^2(t - a_k)_+^2  \Bigr)\qquad\qquad  (0\le t \le a_1+a_2).
$$

\begin{figure}[h]
\begin{center}
\includegraphics [width=0.40 \textwidth] {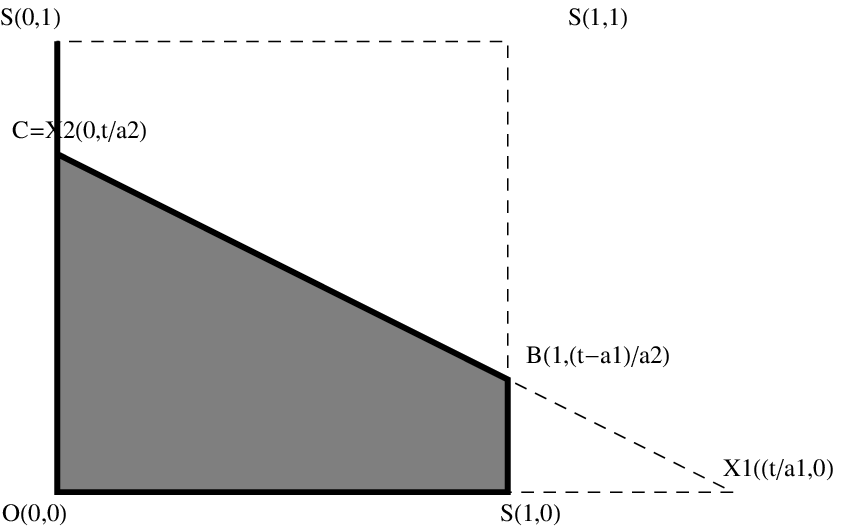}\qquad
\includegraphics [width=0.5 \textwidth]  {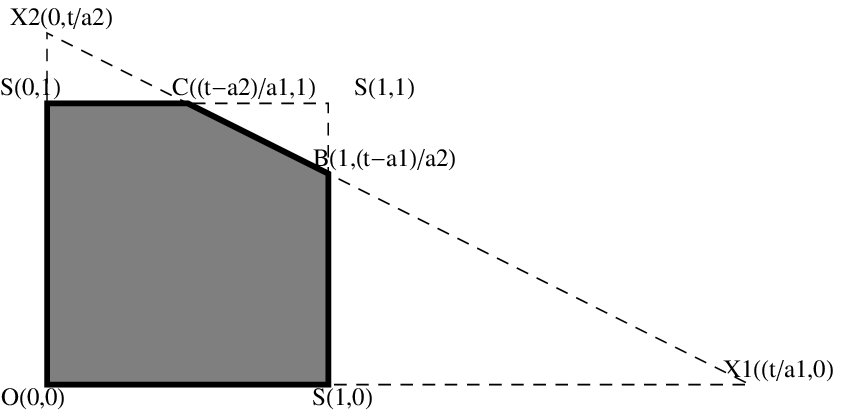}
\caption{Two dimensional case:\ a)  $1\le t/a_1\ \land \ t/a_2\le 1$ \qquad b)  $1\le t/a_1,t/a_2$}
\end{center}
\end{figure}

\noindent(2.4) For $a_1+a_2<t$, it is valid
$$
a_1x_1+a_2x_2\le a_1+a_2<t\qquad (\forall x_1,x_2\in[0,1])
$$
Hence $P_2 = [0,1]^2$, and $p_2=1$.

\smallskip

For arbitrary dimensional space, the hyper-plane $a_1x_1+\cdots+a_mx_m=t$ passes through the points
$$
X_1(t/a_1,0,0,\ldots,0),\  X_2(0,t/a_2,0,\ldots,0),\ldots ,X_m(0,0,\ldots,0,t/a_m).
$$
For $0\le t\le a_1$, the measure of the $m$--dimensional solid between this hyper-plane and coordinate hyper-planes is
$$
\mu(X_m\ldots X_1O)
=\int_0^{\frac{t}{a_1}}dx_1\int_0^{\frac{t-a_1x_1}{a_2}}dx_2\cdots\int_0^{\frac{t-a_1x_1-\cdots-a_{m-1}x_{m-1}}{a_m}}dx_m,
$$
i.e.,
$$
\mu(X_m\ldots X_1O) = \frac{t^m}{m!\prod_{j=1}^m a_j}\,.
$$
Applying the formula (\ref{formula2}), we get
\begin{equation}\label{f0}
f_{0}(t;a_1)  = \sum_{m=0}^\infty \frac{t^m}{m!\prod_{j=1}^m a_j}\qquad (t\in[0,a_1]).
\end{equation}

\section{The analytical approach}

Let us denote by
\begin{equation}\label{expect1}
f_{n}(t;a_k)  = \sum_{i=n}^\infty \frac{t^i}{i!\prod_{j=k}^{i+k-1} a_j}\qquad (t\in[0,a_k])(n\in\mathbb N_0;k\in\mathbb N).
\end{equation}
We accept $f_{n}(t;a_k)\equiv 0\ (t<0)$.

It is important to note that the functions $f_{n}(t;a_k)$ are defined for every $t>0$ according to D'Alambert criteria and the the fact $\inf \{a_i\}=a_1>0$.

\begin{lem}
The next relations are true for all $k,n\in\mathbb N$:
\begin{equation}\label{expect111}
\aligned
f_{n}(t;a_k)   &= \frac{t^n}{n!\prod_{i=k}^{n+k-1} a_{i}} + f_{n+1}(t;a_{k})
,\\
f^{\prime}_{0}(t;a_k)   &= \frac{1}{a_{k}} f_{0}(t;a_{k+1}),\\
f^{\prime}_{n}(t;a_k)   &= \frac{1}{a_{k}} f_{n-1}(t;a_{k+1})\\
f_{0}(t;a_{k+n}) &=     f^{(n)}_0(t;a_k) \prod_{i=k}^{k+n-1} a_i \\
f_0(t;a_{n+1})  &= f^{(n)}_n(t;a_1) \prod_{i=1}^n a_i.
\endaligned
\end{equation}
\end{lem}

The function $f(t,a)$ is known on $(0,a_1)$. It is
$$
f(t,a) = f_0(t;a_1).
$$

\begin{lem}\label{lem32}
It is valid
\begin{equation*}
f(t,a) =1+ f_{1,+}(t;a_1)-f_{1,+}(t-a_1;a_1) \qquad (a_1<t\le a_2).
\end{equation*}
\end{lem}
{\it Proof.} If we take $a_1<t\le a_2$, and a number $x_1\in[0,1]$, then
$$
f(t-a_1x_1,a) = E\{m\in\mathbb N: S_{2,m-1}\le t-a_1x_1 <  S_{2,m} \}=f_{0,+}(t-a_1x_1;a_2).
$$
Including the average of  $f_{0,+}(t-a_1x_1;a_2)$ when $x_1\in[0,1]$, we have
\begin{equation}\label{expect2}
f(t,a) = 1 +  \int_0^1 f_{0,+}(t-a_1x_1;a_2)\ dx_1\qquad (a_1<t\le a_2).
\end{equation}
By change $z=t-a_1x_1$, and the fact
$
f_{0,+}(z;a_2) = a_1 f^{\prime}_{1,+}(z;a_1),
$
yield
\begin{equation*}
f(t,a) = 1 + \frac1{a_1} \int_{t-a_1}^t f_{0,+}(z;a_2)\ dz = 1 +  \int_{t-a_1}^t f^{\prime}_{1,+}(z;a_1)\ dz,
\end{equation*}
what finishes the proof.$\Box$

\begin{lem}\label{lem33} If $a_2<t\le a_3$, then
\begin{equation*}
f(t,a) = 2 + f_{2,+}(t;a_1)-\sum_{k=1}^2 f_{2,+}(t-a_k;a_1)+f_{2,+}(t-a_1-a_2;a_1).
\end{equation*}
\end{lem}
{\it Proof.} If we take two numbers $x_1, x_2\in[0,1]$, then
$$
f(t-a_1x_1-a_2x_2,a) = E\{m\in\mathbb N: S_{3,m-1}\le t-a_1x_1-a_2x_2 < S_{3,m}\},
$$
i.e.,
$$
f(t-a_1x_1-a_2x_2,a) = f_{0,+}(t-a_1x_1-a_2x_2;a_3).
$$
Including the average of  $f_{0,+}(t-a_1x_1-a_2x_2;a_3)$ when $x_1, x_2\in[0,1]$, we have
\begin{equation}\label{expect2}
f(t,a) = 2 +  \int_0^1 dx_1 \int_0^1 f_{0,+}(t-a_1x_1-a_2x_2;a_3)\ dx_2\qquad (a_1+a_2<t\le a_3).
\end{equation}
By change $z=t-a_1x_1-a_2x_2$, we yield
\begin{equation*}
f(t,a) = 2 + \frac1{a_2}\int_0^1 dx_1 \int_{t-a_1x_1-a_2}^{t-a_1x_1} f_{0,+}(z;a_3)dz.
\end{equation*}
Since
$$
f_{0,+}(z;a_3) = a_2 f^{\prime}_{1,+}(z;a_2) = a_1a_2 f^{\prime\prime}_{2,+}(z;a_1),
$$
we have
\begin{equation*}
f(t,a) = 2 +  \int_0^1 f^{\prime}_{2,+}(t-a_1x_1;a_1)dx_1 -   \int_0^1 f^{\prime}_{2,+}(t-a_2 - a_1x_1;a_1)dx_1.
\end{equation*}
Applying the Lemma~\ref{lem32}, we finish the proof. $\Box$

Notice that the conclusions from the previous lemmas can be written in the following forms:
\begin{equation*}
f(t,a) = 1 +  \sum_{k=0}^1 (-1)^k\sum_{0\le i_0<l_1\le 1}
f_{1,+}\Bigl(t-\sum_{j=1}^k a_{i_j};a_1\Bigr)\qquad (a_1<t\le a_2)
\end{equation*}
\begin{equation*}
f(t,a) = 2 +  \sum_{k=0}^2 (-1)^k\sum_{0\le i_0<l_1<i_2\le 2}
f_{2,+}\Bigl(t-\sum_{j=1}^k a_{i_j};a_1\Bigr)\qquad (a_2<t\le a_3)
\end{equation*}

\begin{thm}
 If $a_n<t<a_{n+1}$, it is valid
\begin{equation}\label{minterval}
f(t,a) = n +  \sum_{k=0}^n (-1)^k\sum_{0\le i_0<\cdots <i_k\le k}
f_{n,+}\Bigl(t-\sum_{j=1}^k a_{i_j};a_1\Bigr)
\end{equation}
\end{thm}
{\it Proof.} We will prove by mathematical induction. For $n=1$ it is proven in Lemma~\ref{lem32}.
Also, for $n=2$ it is proven in Lemma~\ref{lem33}.

Let us suppose that the formula (\ref{minterval}) is true for $n$.
If we take numbers $x_1, x_2,\ldots, x_{n+1}\in[0,1]$, then
$$
f(t-S_{1,n+1},a) = E\{m\in\mathbb N: S_{n+2,m-1}\le t-S_{1,n+1} < S_{n+2,m}\},
$$
i.e.,
$$
f(t-S_{1,n+1},a) = f_{0,+}(t-S_{1,n+1};a_{n+2})  \qquad (n+1<m).
$$
Including the average of $f_{0,+}(t-S_{1,n+1};a_{n+2})$  when $x_1, x_2,\ldots, x_{n+1}\in[0,1]$, we have
\begin{equation*}
f(t,a) = n+1 +  \int_0^1 dx_1\cdots \int_0^1 dx_n \int_0^1 f_{0,+}(t-a_1x_1-\cdots -a_{n+1}x_{n+1};a_{n+2})dx_{n+1},
\end{equation*}
By change $z=t-a_1x_1-\cdots -a_{n+1}x_{n+1}$, we yield
\begin{equation*}
\aligned
f(t,a) &= n+1+  \\
&\frac{1}{a_{n+1}}\int_0^1 dx_1 \cdots \int_0^1 dx_n \int_{t-a_1x_1-\cdots-a_{n}x_{n}-a_{n+1}}^{t-a_1x_1-\cdots-a_{n}x_{n}} f_{0,+}(z;a_{n+2})\ dz ,
\endaligned
\end{equation*}
Applying the relations (\ref{expect111}), we have
$$
f_{0,+}(z;a_{n+2})  = f^{(n+1)}_{n+1,+}(z;a_1) \prod_{i=1}^{n+1} a_i,
$$
wherefrom
$$
\aligned
f&(t,a) = n+1\\
&+ \prod_{i=1}^{n}a_i \Biggl(
\int_0^1 dx_1 \cdots \int_0^1 dx_{n-1}\int_0^1 f^{(n)}_{n+1,+}(t-a_1x_1-\cdots-a_{n}x_{n};a_{1}) dx_{n} \\
&- \int_0^1 dx_1 \cdots \int_0^1 dx_{n-1}\int_0^1 f^{(n)}_{n+1,+}(t-a_{n+1}-a_1x_1-\cdots-a_{n}x_{n};a_{1}) dx_{n}\Biggr).
\endaligned
$$
The integrals are average values for the functions  $f_{n+1,+}(t-a_1x_1-\cdots-a_{n}x_{n};a_{1})$ and  $f_{n+1,+}(t-a_{n+1}-a_1x_1-\cdots-a_{n}x_{n};a_{1})$
when $x_1, x_2,\ldots, x_{n}\in[0,1]$.That is why we can apply inductional assumption. Hence
$$
\aligned
f(t,a) &= n+1
+ \Biggl(n-1 + \sum_{k=0}^{n-1} (-1)^k\sum_{0\le i_0<\cdots <i_k\le k} f_{n+1,+}\Bigl(t-\sum_{j=1}^k a_{i_j};a_1\Bigr)\Biggr) \\
&- \Biggl(n-1 + \sum_{k=0}^{n-1} (-1)^k\sum_{0\le i_0<\cdots <i_k\le k} f_{n+1,+}\Bigl(t-a_{n+1}-\sum_{j=1}^k a_{i_j};a_1\Bigr)\Biggr).
\endaligned
$$
By simplifying, we get the formula (\ref{minterval}) for $n+1$, what finishes the proof. $\Box$

\section{Examples}

We prepared the programs in the software {\it Mathematica} to test our conclusions. We used the function
{\tt RandomReal[{0,1},WorkingPrecision->32]}. Blue points  are provided on the next figures in this way every after $500$ random chosen numbers.

\subsection{The Laguerre-type exponential function}

The {\it Laguerre-type exponential function} were introduced by G. Dattoli \cite{Dattoli}:
$$
e_s(t) = \sum_{n=0}^\infty \frac{t^n}{(n!)^{s+1}} \quad (0\le t\le 1)(s\in\mathbb N_0).
$$
From the previous consideration, we cane establish the next statement.
\begin{cor}
When $a_k=k^s\ (k\in\mathbb N)(s\ge 0)$, it is valid
\begin{equation*}
f(t,a) =
\begin{cases}
e_s(t), & (0\le t\le 1)\\
 e_s(t)-e_s(t-1)+1 , & (1<t\le 2^s).
 \end{cases}
\end{equation*}
In general, in the formula (\ref{minterval}) appears the next fuction
\begin{equation*}
f_{n,+}(t) = e_s(t)- \sum_{k=0}^{n-1} \frac{t^i}{(i!)^{s+1}}.
\end{equation*}
\end{cor}

\begin{rem}\rm
The case $s=0$ was examined in \cite{Curgus}. Since $a=\{1\}$,  we easily find that $f(t,a)=e^t$ for $0\le t \le 1$, but our approach can not give $f(t,a)$ for $t>1$.
\end{rem}

\begin{center}
\begin{figure}[h]
\includegraphics [width=0.45 \textwidth] {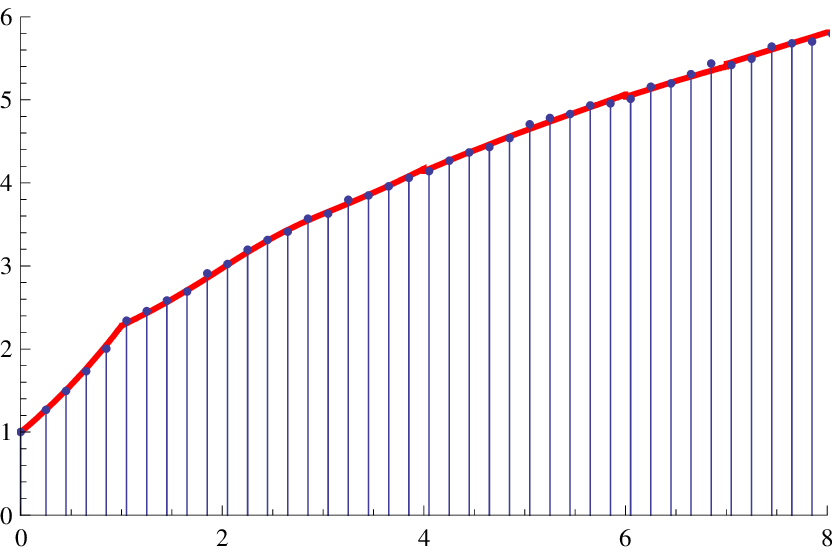}\qquad
\includegraphics [width=0.45 \textwidth] {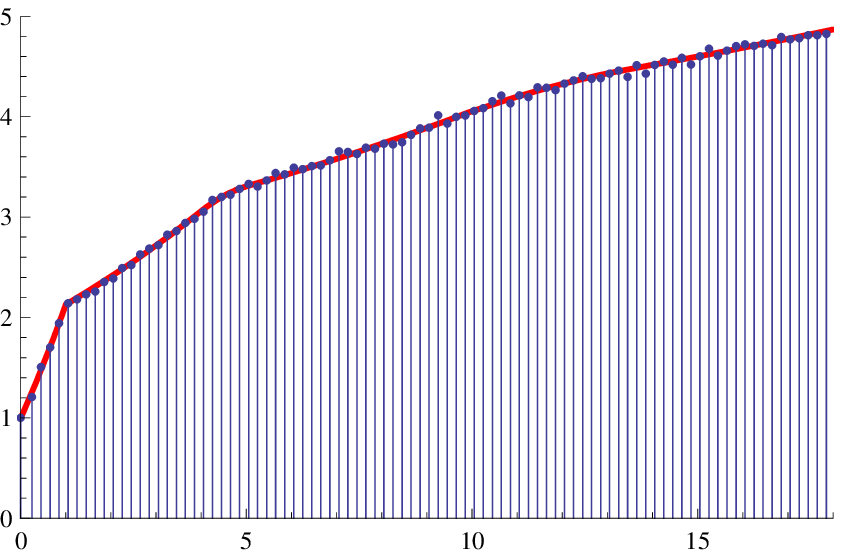}
\caption{Dattoli cases:\ a) $a_k = k$ \qquad\qquad \qquad\qquad \qquad b) $a_k = k^2$.}
\end{figure}
\end{center}

\subsection{The Laguerre-type $q$--exponential function}

Another theoretically interesting example is the case
$$
a_k=1-q^k   \qquad (0\le q<1) \ (k\in\mathbb N).
$$
It is leading to a {\it Laguerre type  $q$--exponential function}\cite{Stankovic}
$$
e_{1,q}(t) = \sum_{n=0}^\infty \frac{t^n}{n!(q;q)_n}\quad (0<t\le 1-q).
$$
Here
$$
(q;q)_0=1,\qquad (q;q)_n=\prod_{k=1}^n(1-q^k) \quad (n\in\mathbb N).
$$
\begin{figure}[h]
\begin{center}
\includegraphics [width=0.5 \textwidth] {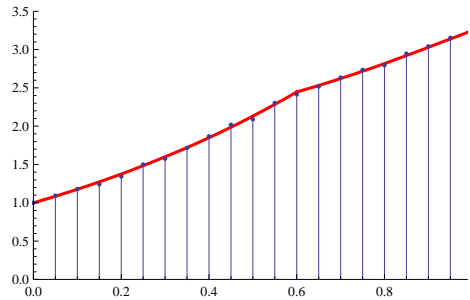}\qquad
\caption{ Case with $q$-numbers $a_k = 1-q^k$ \quad $(q = 0.4)$.}
\end{center}
\end{figure}

\noindent {\bf Acknowledgements.} This paper is supported by the Ministry of Science and Technological Development of the Republic Serbia, projects No $174011$ and No $44006$.

\end{document}